\documentclass{article}
\usepackage{graphicx} 

\usepackage{amsfonts,amssymb,epsfig}
\usepackage{amsmath}
\usepackage{amsthm}
\theoremstyle{plain}
\date{}

\newtheorem*{thm1}{Theorem 1}

\def\e{\varepsilon}

\title{Entropy versus influence for complex functions of modulus one 2
}

\author{Gideon Schechtman } 
\begin{document}
\maketitle

 \centerline{\date{April 2024}}

\begin{abstract}
This is a simplification of a previous version of this ArXiv note. We present an example of a function $f$ from $\{-1,1\}^n$ to the unit sphere in $\mathbb{C}$ with influence bounded by $1$ and entropy of $|\hat f|^2$ larger than $\frac12\log n$. 

\end{abstract}

\section{Introduction}
We denote by $\e_i:\{-1,1\}^n\to \{-1,1\}$ the projection onto the $i$-s coordinate: $\e_i(\delta_1,\dots,\delta_n)=\delta_i$. For a subset $A$ of $[n]:=\{1,\dots,n\}$ we denote $W_A=\prod_{i\in A}\e_i$, $W_A:\{-1,1\}^n\to \{-1,1\}$. The $W_A$-s are the characters of the Cantor group $\{-1,1\}^n$ (with coordintewise multiplication) and form an orthonormal basis in $L_2$ of the Cantor group equipped with the normalized counting measure. In this note we shall be concerned with functions from $\{-1,1\}^n$ into the complex plane, $\mathbb{C}$. These can also be considered as a couple of real functions. Each such function $f: \{-1,1\}^n\to \mathbb{C}$ has a unique expansion
\[
f=\sum_{A\subseteq [n]}\hat f(A)W_A,
\]
 where $\hat f(A)\in \mathbb{C}$ are given by
 \[
 \hat f(A)=\mathbf{E}(fW_A)=2^{-n}\sum_{\e\in\{-1,1\}^n}f(\e)W_A(\e).
 \]
 Note that if $f: \{-1,1\}^n\to \mathbb{R}$, then $\hat f(A)\in \mathbb{R}$ for every $A\subseteq[n]$. The orthonormality of the $W_A$-s implies that
 \[
 \|f\|_2^2:=2^{-n}\sum_{\e\in\{-1,1\}^n}|f(\e)|^2=\sum_{A\subseteq[n]}|\hat f(A)|^2.
 \]

 Define the influence of a function $f: \{-1,1\}^n\to \mathbb{C}$ by
 \begin{equation}\label{eq:defI}
 I(f)=\sum_{A\subseteq[n]}|\hat f(A)|^2|A|
 \end{equation}
 where for $A\subseteq[n]$, $|A|$ denotes the cardinality of $A$. This object, especially for boolean functions, is a deeply studied one and quite influential (but this is not the reason for the name...) in several directions. We refer to \cite{odonnell} for some information. A recent paper dealing with
 the subject is  \cite{kklms}.

 The entropy of the sequence $|\hat f(A)|^2$ is given by
 \begin{equation}\label{eq:defH}
 H(|\hat f(A)|^2)=-\sum_{A\subseteq[n]}|\hat f(A)|^2\log |\hat f(A)|^2.
 \end{equation}
 ($0\log 0:=0$). The base of the $\log$ does not really matter here. For concreteness we take the $\log$ to base $2$. Note that if $f$ has $L_2$ norm $1$ then the sequence  $\{|\hat f(A)|^2\}_{A\subseteq[n]}$ sums up to $1$ and thus this is the usual definition of entropy of this probability distribution.

 The entropy influence conjecture of Friedgut and Kalai \cite{fk} is that
 for some absolute constant $K$, for all $n$ and all boolean functions $f:\{-1,1\}^n\to\{-1,1\}$
 \[
 H(|\hat f(A)|^2)\le KI(f).
 \]
 For the significance of this conjecture we refer to the original paper \cite{fk}, and to Kalai's blog \cite{k} (embedded in Tao's blog) which reports on all significant results concerning the conjecture. \cite{kklms} establishes a weaker version of the conjecture. Its introduction is also a good source of information on the problem.

In version 1 of this note, which can still be found on the ArXiv, we showed that the analogous version of the conjecture for complex functions on  $\{-1,1\}^n$ which have modulus $1$ fails. This solves a question raised by Gady Kozma some time ago (see \cite{k}, comment from April 2, 2011). More specifically, we proved

\begin{thm1}
For each $n=1,2,\dots$ there is a function $f:\{-1,1\}^n\to \{z\in \mathbb{C}; \ |z|=1\}$ with
\[
I(f)<1, \ \ {\mbox {and}}\ \ H(|\hat f|^2)>\frac{n}{n+1}\log n.
\]
\end{thm1}

Here we give an embarrassingly simple presentation of an example of such a function (although it can be shown to be a version of the example in the previous version of this note). As was written in the previous version, an anonymous referee of version 1 wrote that the theorem  was known to experts but not published. Maybe the presentation below is what was known. \\
Consequently, this note will not be sent for publication.

Let $n\in{\mathbb N}$  and consider
\[F=(1+1/n)^{-n/2}\prod_{j=1}^n(1+\imath \epsilon_j/\sqrt n)\]
Clearly,
\[|F|\equiv 1.\]
The Walsh--Fourier expansion of F is
\[
F=(1+1/n)^{-n/2}\sum_{A\subseteq\{1,\dots,n\}}\left(\frac{\imath}{\sqrt n}\right)^{|A|}W_A.
\]
consequently,
\[
I(F)=(1+1/n)^{-n}\sum_{k=1}^n kn^{-k}{n\choose k}=(1+1/n)^{-n}(1+1/n)^{n-1}=n/(n+1).
\]
and
\begin{multline*}
H(F)=(1+1/n)^{-n}\sum_{k=1}^n \log(n^{k}(1+1/n)^n)n^{-k}{n \choose k}\\
>(1+1/n)^{-n}\log n \sum_{k=1}^n k n^{-k}{n \choose k}=(n\log n)/(n+1).
\end{multline*}

%
%

\noindent G. Schechtman\\
Department of Mathematics\\
Weizmann Institute of Science\\
Rehovot, Israel\\
{\tt gideon@weizmann.ac.il}

\end{document}